\documentclass[12pt]{article}

\usepackage[colorlinks=true, pdfstartview=FitV, linkcolor=blue, citecolor=blue, urlcolor=blue]{hyperref}

\usepackage{amssymb,amsmath, amscd}
\usepackage{times, verbatim, url}
\usepackage{graphicx}

\DeclareFontFamily{OT1}{rsfs}{}
\DeclareFontShape{OT1}{rsfs}{n}{it}{<-> rsfs10}{}
\DeclareMathAlphabet{\mathscr}{OT1}{rsfs}{n}{it}

\newtheorem{theorem}{Theorem}[section]

\DeclareMathOperator{\Gal}{Gal}

\DeclareMathOperator{\SL}{SL}

\DeclareMathOperator{\Tr}{Tr}
\DeclareMathOperator{\disc}{disc}

\DeclareMathOperator{\Sym}{Sym}

\DeclareMathOperator{\ord}{ord}

\DeclareMathOperator{\N}{N}

\def\R{{\mathbb R}}
\def\F{{\mathbb F}}
\def\PP{{\mathbb P}}
\def\bC{{\mathbb C}}
\def\Q{{\mathbb Q}}
\def\Z{{\mathbb Z}}
\def\A{{\mathbb A}}
\def\O{{\mathscr O}}
\def\p{{\frak p}}
\def\q{{\frak q}}

\def\sL{{\mathscr L}}

\title{Two encounters with the $p$--adic Stark conjecture}

\author{Benedict H Gross and Samit Dasgupta}

\begin{document}
\maketitle

\tableofcontents

\newpage

In this note we describe our personal encounters with the $p$-adic Stark conjecture. Gross describes the period between $1977$ and $1986$ when he came to formulate these conjectures (\S\ref{s:gs}--\ref{s:integral}), and Dasgupta describes the period between $1998$ and $2021$, when he worked with others to finally prove them (\S\ref{s:exact}--\ref{s:ic}).

\section{Gauss sums and Stickelberger's theorem} \label{s:gs}

In my third year of graduate school, while casting about for a thesis, I (Benedict Gross) spent a lot of time talking to David Rohrlich and Neal Koblitz. who had just arrived as assistant professors. We were trying to see if some of the methods in Barry Mazur's great paper  ``Modular Curves and the Eisenstein Ideal" could be transferred to a study of the Fermat curve.
David had determined the analog of the cuspidal group in his PhD thesis \cite{R}, which was a start.
We soon realized that we needed to understand the $L$-function, which involved the theory of Gauss and Jacobi sums \cite{GR},\cite{W}.\\

Recall Euler's congruence for the quadratic residue symbol at an odd prime $p$: 
$$\big(\frac{x}{p}\big) \equiv x^{(p-1)/2} ~~ ({\rm mod}~ p).$$
This generalizes to a congruence for the $N$ power residue symbol, with $N > 2$, as follows. Let $\p$ be a prime ideal in the ring $A = \Z[\mu_N]$ and let $q = p^f$ be the number of elements in the residue field $\F_{\p} = A/\p$. We will always assume that the residual characteristic is prime to $2N$. Then the $N^{th}$ roots of unity in $A$ map isomorphically to the $N^{th}$ roots of unity in $\F_{\p}$ under reduction modulo $\p$. Hence $q \equiv 1$ modulo $N$.

For any non-zero element $x$ in the residue field, the power $x^{(q-1)/N}$ has order $N$ in the (cyclic) group $\F_{\p}^*$. Hence there is a unique $N^{th}$ root of unity in the ring $A$, denoted $\big(\frac{x}{\p}\big)_N$,  which satisfies
$$\big(\frac{x}{\p}\big)_N \equiv  x^{(q-1)/N} ~~ ({\rm mod}~ \p).$$
The map $\chi(x) = \big(\frac{x}{\p}\big)_N$ is a multiplicative character of $\F_{\p}^*$ taking values in $\mu_N \subset A^*$.\\

The Gauss sums in question have the form
$$g(a,\p) = -\sum_{x \in \F_{\p}^*} \chi(x)^{-a}~~ \psi(\Tr(x))$$
where $a$ is a non-zero class in $\Z/N\Z$ and $\psi: \F_p \rightarrow \mu_p$ is a non-trivial additive character. Since each term in the sum is the product of an $N^{th}$ root of unity and a $p^{th}$ root of unity, the sum $g(a,\p)$ is an algebraic integer in the number field $\Q(\mu_N,\mu_p)$.
A simple argument using Galois theory shows that the Jacobi sum
$$J(a,\p) = g(a,\p)^N$$
is an element of the ring $A = \Z[\mu_N]$ of integers in $\Q(\mu_N)$, which is independent of the choice of additive character $\psi$.\\

Stickelberger's theorem gives the factorization of the principal ideal of $A$ generated by $J(a,\p)$.  Since 
$$g(a,\p)\overline{g(a,\p)} = q = p^f$$
this factorization only involves primes dividing the rational prime $p$.  Let $\zeta$ be a primitive $N^{th}$ root of unity and for $u \in (\Z/N\Z)^*$ let $\sigma_u$ be the automorphism of $A = \Z[\zeta]$ which maps $\zeta$ to $\zeta^u$. Since $\sigma_p$ generates the decomposition subgroup of the prime $\p$, the primes dividing $p$ in $A$ have the form $\sigma_u(\p)$, with $u \in (\Z/N\Z)^*/\langle p \rangle$.

For a non-zero class $a$ in $\Z/N\Z$ we let $\langle a \rangle$ be the unique integer between $0$ and $N$ which is congruent to $a$ modulo $N$. Then \cite{W}:
$$(J(a,\p)) = \prod_{u \in (\Z/N\Z)^*} \sigma_u^{-1}(\p)^{\langle ua \rangle}.$$
In particular:
$\ord_{\p} J(a,\p) = \langle a \rangle + \langle pa \rangle + \ldots + \langle p^{f-1}a \rangle.$

Since $J(a,\p) = g(a,\p)^N$, this follows from a similar result on Gauss sums. We view $g(a,\p)$ in the completion of $\Q(\mu_N, \mu_p)$ at the unique prime dividing $\p$. Then $g(a,\p)$ lies in the subfield $\Q_p(\mu_p)$. This is a Kummer extension of $\Q_p$, obtained by adjoining the $(p-1)^{st}$ roots of $-p$. Let $\zeta_p = \psi(1)$ and let $\pi$ be the unique $(p-1)^{st}$ root of $-p$ which is congruent to $(\zeta_p -1)$ modulo $(\zeta_p -1)^2$. Then
$$g(a,\p) = \pi^{(p-1)(\langle a \rangle + \langle pa \rangle+ \ldots + \langle p^{f-1}a \rangle)/N}~ E(a,\p),$$
where the exponent of $\pi$ is an integer and $E(a,\p)$ is a unit in $\Z_p^*$.

Stickelberger found a formula
for this unit modulo $p$. Since $\langle a \rangle/N$ is a rational number between $0$ and $1$, the integer $(p^f - 1) \langle a \rangle/N$ has a $p$-adic expansion
$$(p^f - 1) \langle a \rangle/N = z_0 + z_1p + z_2p^2 + \ldots z_{f-1}p^{f-1}$$
with $0 \leq z_i < p$. Then (cf. \cite[Ch 1]{L})

$$\ord_{\pi}(g(a,\p)) = z_0 + z_1 + \ldots + z_{f-1},$$
$$E(a,\p) \equiv 1/z_0!z_1!\ldots z_{f-1}! ~~({\rm mod} ~p).$$

\section{The $p$-adic Gamma function and $p$-adic $L$-functions}

Neal  Koblitz and I found a extension of Stickelberger's factorial formula for $E(a,\p)$, using the $p$-adic Gamma function \cite{GK}. Our argument relied crucially on a limit formula for Jacobi sums, which Nick Katz had obtained from a study of the crystalline cohomology of the Fermat curve (cf. \cite[Ch 15]{L}).

The function $\Gamma_p(x)$
 is a continuous function from $\Z_p$ to $\Z_p^*$, which is determined by its values on (the dense set) of positive integers $m$. For $1 \leq m \leq p$ we have $\Gamma_p(m) = (-1)^m(m-1)!$ similar to the classical Gamma function. In general, $\Gamma_p(m)$ is $(-1)^m$ times the product of the positive integers less than $m$ which are prime to $p$. The function $\Gamma_p$ is continuous in the $p$-adic topology (by Wilson's theorem) and the values of $\Gamma_p$ at negative integers can be obtained from the formula
 $$\Gamma_p(z)\Gamma_p(1-z) = (-1)^{m(z)}$$
 where $m(z)$ is the unique integer with $0 < m(z) \leq p$ which is congruent to $z$ modulo $p$.
  
Our formula for the unit $E(a,\p)$ involves the values of $\Gamma_p$ at rational numbers between $0$ and $1$ with denominator $N$. Recall that the residue field $A/\p$ has cardinality $q = p^f$. Then
 $$E(a,\p) = \prod_{n=0}^{f-1} \Gamma_p(\langle p^n a \rangle/N).$$
 
 At about the same time that we found this formula, Ralph Greenberg and Bruce Ferrero found a formula relating the $p$-adic logarithms of the {\bf same} values of the $p$-adic Gamma function to the first derivatives of $p$-adic $L$-functions \cite{FG}. These $p$-adic $L$-functions were introduced by Kubota and Leopoldt, and are associated to Dirichlet characters. Let $\chi: (\Z/N\Z)^* \rightarrow \bC^*$ be an odd Dirichlet character of conductor $N$. Then the complex Dirichlet $L$-function $L(\chi,s)$ does not vanish at $s = 0$ and 
 $$L(\chi,0) = - \sum_{a = 1}^N \chi(a) (a/N) = -B_{1,\chi},$$
 $$L'(\chi, 0) = \sum_{a = 1}^N \chi(a) \log \Gamma(a/N) + B_{1,\chi} \log N.$$
 
Greenberg and Ferrero found an analog of these formulas for the $p$-adic L-function associated to the even $p$-adic character $\chi\omega$, where $\omega$ is the Teichmuller character:
$$L_p(\chi \omega,0) = -(1 - \chi(p)) B_{1,\chi},$$
$$L'_p(\chi \omega, 0) = \sum_{a = 1}^N \chi(a) \log_p \Gamma_p(a/N) + (1 - \chi(p))B_{1,\chi} \log_p(N),$$
where $\log_p$ is the Iwasawa logarithm, with $\log_p(p) = 0$.
Since $B_{1,\chi} \neq 0$, the value at $s = 0$ is non-zero precisely when $\chi(p) \neq 1$.

 When $\chi(p) = 1$, the derivative is given by the simpler formula
$$L'_p(\chi \omega, 0) = \sum_{a = 1}^N \chi(a) \log_p \Gamma_p(a/N) = \sum_{a \in (\Z/N\Z)^*/<p>} \chi(a) \log_p \prod_{n=0}^{f-1} \Gamma_p(\langle p^n a \rangle/N).$$
Using the relation between the $p$-adic Gamma function and Gauss sums, and the fact that the Iwasawa logarithm satisfies $\log_p(p) = 0$, this becomes
$$L'_p(\chi\omega,0) =  \sum_{a \in (\Z/N\Z)^*/<p>} \chi(a) \log_p g(a,\p) = (1/N)\sum_{a \in (\Z/N\Z)^*/<p>} \chi(a) \log_p J(a,\p).$$

Greenberg and Ferrero then use an argument in linear forms in logarithms of algebraic numbers to prove that $L'_p(\chi\omega,0)  \neq 0$. (We will do a simple case in the next section.) Hence the order of $L'_p(\chi\omega,s)$ at $s = 0$ is zero when $\chi(p) \neq 1$ and is one when $\chi(p) = 1$. This had been predicted by an earlier paper by Greenberg on Iwasawa modules \cite{G}.

\section{The imaginary quadratic case}

Consider the case when the odd Dirichlet character $\chi$ is quadratic, of conductor $D$. Then $\chi$ determines the splitting of primes in the imaginary quadratic field $K = \Q(\sqrt {-D})$. Indeed, a prime $p$ splits if and only if $-D$ is a (non-zero) square modulo $p$. But then quadratic reciprocity implies that $p$ splits if and only if $\chi(p) = 1$.

From a consideration of Euler products, we obtain a factorization 
$$\zeta_K(s) = \zeta_{\Q}(s) L(\chi,s).$$
But $\zeta_{\Q}(0) = -1/2$ and $\zeta_K(0) = -h/w$, where $h$ is the class number and $w$ is the order of the unit group of the ring of integers $A_K$. (These identities are a consequence of a more general formula describing the zeta function at $s = 0$, which we will review in the next section.) Using these two values, we obtain Dirichlet's beautiful formula for the class number of $K$:
$$L(\chi,0) = 2h/w = -B_{1,\chi} = - \sum_{a = 1}^D \chi(a) (a/D) =  \sum_{\chi(a) = -1} \langle a \rangle /D  ~~ - \sum_{\chi(a) = +1} \langle a \rangle /D.$$
where the last two sums are taken over $a \in(\Z/D\Z)^*$.

In the $p$-adic case, when $\chi(p) = -1$, the value at $s = 0$ is
$L_p(\chi \omega,0) = 4h/w.$
On the other hand, when $\chi(p) = +1$, the $p$-adic L-function vanishes at $s = 0$ and
$$L'_p(\chi \omega,0) = (1/D)\sum_{a \in (\Z/D\Z)^*/<p>} \chi(a) \log_p J(a,\p) = (1/D)\log_p \Big( \prod_{\chi(a) = +1} J(a,\p)/ \prod_{\chi(a) = -1} J(a,\p)\Big),$$
where the products are now taken over $a$ in the quotient group $(\Z/D\Z)^*/<p>$.

The imaginary quadratic field $K = \Q(\sqrt{-D})$ is a subfield of $\Q(\mu_D)$, and the Galois group of $\Q(\mu_D)/K$ is the subgroup of $a \in (\Z/D\Z)^*$ with $\chi(a) = +1$. Hence the product of Jacobi sums in the expression for the derivative lies in $K$. This product generates a power of the ideal $\frak q/\overline{\frak q}$, where $\frak q$ is the factor of $p$ in $K$ which is divisible by $\p$.
We can determine which power by using the ideal factorization of the Jacobi sums and the class number formula for $B_{1,\chi}$:
$$\ord_{\p} \Big( \prod_{\chi(a) = +1} J(a,\p)/ \prod_{\chi(a) = -1} J(a,\p)\Big) = \sum_{\chi(a) = +1}\langle a \rangle - \sum_{\chi(a) = -1} \langle a \rangle = -D2h/w.$$
If we let $\alpha$ generate the principal ideal $(\frak q)^h$, then
$$L'_p(\chi \omega,0) = (-2/w) \log_p(\alpha/\overline{\alpha}) = (4/w) \log_p(\overline{\alpha}),$$
where $K$ is embedded in $\Q_p$ via the place $\frak q$. Since $\overline{\alpha}$ is a unit in this embedding, which is not a root of unity, $\log_p( \overline{\alpha} )\neq 0$ and the first derivative is non-zero.

\section{Stark's conjecture for Artin L-series}

At the same time that I was learning about the results of Ferrero and Greenberg, I was attending a series of lectures at Harvard on Stark's conjectures \cite{S}.  Harold Stark came over from MIT to give three introductory lectures, which included his recent computations involving the L-functions of holomorphic modular forms of weight $1$. Then John Tate presented the conjectures from another point of view \cite{T}.

Stark's conjectures concern the behavior of Artin L-functions at $s = 0$, and seek to generalize the work of Dirichlet and Dedekind on zeta functions. Let $K$ be a number field with $r_1$ real embeddings and $r_2$ inequivalent complex embeddings. Then the order of $\zeta_K(s)$ at $s = 0$ is equal to $r = (r_1 + r_2 -1)$. This is also equal to the rank of the unit group $A^*$ of the ring of integers $A$. Furthermore, the leading term in the Taylor expansion of $\zeta_K(s)$ at $s = 0$ is equal to $(-hR/w) . s^r$, where $h$ is the class number of $A$, $w$ is the number of roots of unity in $A^*$, and the regulator $R$ is the absolute value of the determinant of an $r \times r$ matrix whose entries are logarithms of units. (When $K$ is imaginary quadratic, as in the previous section, $\zeta_K(0) = -h/w$, as the group $A^*$ is finite and $R = 1$.)

Stark's amazing idea was an attempt to factor the leading term at $s = 0$ in the same way that the zeta function of $K$ factors as a product of Artin L-functions. If $K$ is a Galois extension of $k$ with group $G$, then Artin obtained the factorization
$$\zeta_K(s) = \prod_V L(V,s)^{\dim V}$$
where the product is taken over the irreducible complex representations $V$ of $G$. The function $L(V,s)$ is an Artin L-series, defined in the half plane where the real part of $s$ is greater than one by an Euler product over the finite primes of the field $k$:
$$L(V,s) = \prod_{\frak p} \det(1 - Fr_{\frak p}(N\frak p)^{-s}|V^{I_{\frak p}})^{-1}.$$

To determine the leading term of $L(V,s)$ at $s = 0$ the first thing is to determine the order of vanishing there. For each Archimedian place $v$ of $k$, let $G_v$ be a decomposition subgroup (of order $1$ or $2$) in $G$. If $V$ is the trivial representation, $L(V,s) = \zeta_k(s)$ has a simple pole at $s = 1$. If $V$ is a non-trivial irreducible representation, the function $L(V,s)$ is regular and {\bf non-zero} at $s = 1$. From the Archimedian factors in the functional equation relating $s$ to $1-s$, one concludes that for any $V$:
$$\ord_{s = 0}L(V,s) = \sum_v \dim V^{G_v}  - \dim V^G$$
Call this order of vanishing $r(V)$. 

Let $Y$ be the free abelian group on the Archimedian places $v$ of $K$ and let $X$ be the subgroup of $Y$ consisting of elements $\sum a_v .v$ with $\sum a_v = 0$. As rational representations of $G$, $Y \otimes \Q$ is a sum of permutation representations on the set $G/G_v$ and $X \otimes \Q$ is this direct sum minus the trivial representation. From Frobenius reciprocity, we obtain the formula
$$\dim(X \otimes V)^G = \sum_v \dim V^{G_v}  - \dim V^G = r(V)$$
 Stark conjectured that the primary component of the leading term of $L(V,s)$ at $s = 0$ should be a regulator $R(V)$ -- the determinant of an $r(V) \times r(V)$ matrix of logarithms of units, in the $V$-isotypic component of the units $A^*$ of $K$. In particular, in the special case when $r(V) = 1$ the first derivative $L'(V,0)$ should determine (the logarithm of) a unit in the Galois extension $K$. That opened a completely new approach to Hilbert's  $12^{th}$ problem on the construction of class fields! 

Here is how Tate formulated Stark's main conjecture. Define the $G$ homomorphism
$$\lambda: A^* \rightarrow X \otimes \R$$
by the formula $\lambda(\epsilon) = \sum_v \log |\epsilon|_v .v$ where $| |_v$ is the canonical absolute value at the place $v$. Dirichlet's unit theorem states that the kernel of $\lambda$ is the torsion subgroup (of roots of unity), and that the image is a cocompact subgroup of $X\otimes \R$. The classical regulator for $K$ is the covolume of this subgroup. Another way to state Dirichlet's theorem is that the induced map
$\lambda: A^* \otimes \R \rightarrow X \otimes \R$
is an isomorphism of $\R[G]$-modules. 

Since the two rational $G$ modules $A^* \otimes \Q$ and $X \otimes \Q$ are isomorphic over $\R$, they have the same character and are isomorphic over $\Q$. 
Choose such a $G$-isomorphism
$\phi: X \otimes \Q \rightarrow A^* \otimes \Q$
and for an irreducible complex representation $V$ of $G$ consider the isomorphism of complex vector spaces
$$(1 \otimes \lambda \phi)_V:  (V \otimes X)^G \rightarrow (V \otimes X)^G.$$
The Stark regulator $R(V)$ is the determinant of the map $(1 \otimes \lambda \phi)_V$, and the main conjecture is that the leading term of $L(V,s)$ at $s = 0$ is a product of the form $R(V)A(V)$, where $A(V)$ is an algebraic number in the field generated by the character of $V$.

\section{A conjecture for $p$-adic Artin L-series} \label{s:padicconj}

This was such a beautiful formulation that I began to wonder if it could be transposed to the case of $p$-adic L-functions. In this case, we assume that the base field $k$ is totally real and that the Galois extension $K$ is a CM field, with complex conjugation $c \in G$. We only consider totally odd representations $V$, where $c = -1$. 

If we view $V$ as a complex representation, then $r(V) = 0$ and
$L(V,0) \neq 0$. In fact, by results of Siegel, $L(V,0) = A(V)$ is an algebraic number in the field generated by the character of $V$. 
What is more relevant to the $p$-adic L-function is the partial L-function $L_S(V,s)$, defined as an Euler product over the places not in $S$, where $S$ is the set of places of $k$ dividing $\infty$ and $p$. Hence 
$$L_S(V,s) = \prod_{\frak p} \det(1 - Fr_{\frak p}(N\frak p)^{-s}|V^{I_{\frak p}})~L(V,s) $$
where the finite product is taken over the primes $\frak p$ of $k$ which divide $p$. This function vanishes to order
$$r_p(V) = \sum \dim V^{G_{\frak p}}$$
at $s = 0$ and the leading term has the form
$$\prod_{\frak p} (-\log p)^{\dim V^{G_{\frak p}}}  ~A_p(V)$$
where $A_p(V)$ is an algebraic number in the field generated by the character of $V$, closely related to the value
$L(V,0)$

Henceforth we view $V$ as a representation of $G$ over $\bC_p$. The $p$-adic L-function $L_p(V\omega,s)$ defined by Deligne-Ribet \cite{DR} and Cassou Nogues \cite{C}, conjecturally takes the value
$$L_p(V \omega,0) = \prod_{\frak p} \det(1 - Fr_{\frak p}|V^{I_{\frak p}})~L(V,0) = L_S(V,0)$$
at $s = 0$, where the product is taken over the primes $\frak p$ of $k$ which divide the rational prime $p$. This is known to be true when $G$ is abelian, in which case an odd irreducible representation $V$ is just a character $\chi$ of $G$ with $\chi(c) = -1$ and
$$L_p(\chi\omega,0) = \prod_{\frak p} (1 - \chi(\frak p))~L(\chi,0) = L_S(\chi,0)$$
We note that the proof in \cite{DR} of the congruences necessary to establish the existence of a continuous function $L_p(\chi\omega,s)$ on $\Z_p$ with the correct values on negative integers uses the Fourier expansions of Hilbert Eisenstein series.

In the abelian case, Greenberg conjectured that the order of vanishing of $L_p(\chi\omega)$ at $s = 0$ is equal to
$\#\{\frak p: \chi(\frak p) = 1 \}$ \cite{G}. More generally, I guessed that the order of vanishing of $L_p(V\omega,s)$ at $s = 0$ should be equal to
$r_p(V)$. I also defined a $p$-adic regulator $R_p(V)$, as the determinant of an $r_p(V) \times r_p(V)$ matrix whose entries are the $p$-adic logarithms of ``minus $p$-units" in the extension field $K$, and conjectured that the leading term at $s=0$ is equal to the product $R_p(V)A_p(V)$. That was the earliest version of the $p$-adic Stark conjecture \cite{Gr}, and here is a precise statement.

The minus $p$ units $U^{-1}$ are the $S$-units $a$ of $K$ which satisfy $a^{1+c} = 1$. They have trivial valuation at all infinite places and at the finite primes which do not divide $p$. Let $Y$ be the free abelian group on the places $\frak q$ dividing $p$ in $K$ and let $X$ be the subgroup of those $\sum a(\frak q).\frak q$ in $Y$ with $\sum a(\frak q) = 0$. In analogy to the homomorphism $\lambda$ in the complex case, I defined a $G$-homomorphism
$\lambda_p: U^{-1} \rightarrow \Q_p \otimes X^{-1}$
by the formula
$$\lambda_p (a) = \sum_{\frak q} \log_p (N_{K_{\frak q}/\Q_p}(a)) . ~\frak q$$
and conjectured that this map gives an isomorphism of $G$-modules $\Q_p \otimes U^{-1} \rightarrow \Q_p \otimes X^{-1}$. This is an analog of Leopold's conjecture on units, and is probably just as difficult to establish.

There is an integral $G$-homomorphism
$g: U^{-1} \rightarrow X^{-1}$
which is defined by the formula
$$g(a) = \sum_{\frak q} f_{\frak q}\ord_{\frak q}(a) .~ \frak q$$
where $f_{\frak q}$ denotes the degree of the residue field. This gives an isomorphism of $G$-modules
$\Q \otimes U^{-1} \rightarrow \Q \otimes X^{-1}$, and we define the regulator associated to $V$ by
\begin{equation} \label{e:rpdef}
R_p(V) = \det (1 \otimes \lambda_p g^{-1}: (V \otimes X^{-1})^G).
\end{equation}
Since $X^{-1}$ is the minus eigenspace for $c$ on the permutation representation on primes dividing $p$, the space $(V \otimes X^{-1})^G$ has dimension $r_p(V)$.

In making the general conjecture, I was guided by the special case we considered earlier, where $k = \Q$, $K$ is an imaginary quadratic field, and $V$ corresponds to the non-trivial character $\chi$ of the Galois group. Assume that the rational prime $p$ splits in $K$: $(p) = \frak q \overline{\frak q}$. Then
$L_S(\chi, s) = (1 - \chi(p)p^{-s})L(\chi, s)$ vanishes to order $1$ at $s =0$ and $L_S'(\chi,0) = (\log p) L(\chi,0) = (\log p) (2h/w)$. Hence the algebraic number $A_p(\chi)$ defined above is equal to  $-(2h/w)$. On the other hand, if $\alpha$ is a generator of the principal ideal $(\frak q)^h$, then the ratio $\alpha/\overline{\alpha}$ is a minus $p$ unit in $K$. We have $R_p(\chi) = (1/h) \log_p(\alpha/\overline{\alpha})$, where $\alpha/\overline{\alpha}$ is viewed in the completion $K_{\frak q} = \Q_p$. The conjecture predicts that $L_p'(\chi, 0) = R_p(\chi)A_p(\chi) = (-2/w)\log_p(\alpha/\overline{\alpha})$, which we have seen is true. Using our general results on Gauss sums, I checked that the conjecture was also true whenever $K$ is an abelian extension of $\Q$.

\section{The abelian rank one case} \label{s:abelian}

Stark identified a particularly interesting special case in his conjectures, where the Galois extension $K/k$ is abelian with group $G$ and the $L$-functions $L_S(\chi, s)$ vanish to order one at $s = 0$. This also has a $p$-adic formulation. In this we assume that $k$ is totally real and that $K$ is a CM field. Let $S$ denote a finite set of places of $k$, containing all Archimedean places, all places dividing $p$, and all places ramified in $K$. 

We further assume that there is a place $\frak p$ dividing $p$ which splits completely in $K$, and fix a prime factor $\frak q$ of $\frak p$ in $K$. Since $\p$ is split, we have $r(\chi) = r_p(\chi) \geq 1$ for all characters $\chi$ of $G$. We will consider the first derivatives of both the complex and $p$-adic L functions at $s = 0$. In fact, we will consider the partial zeta functions, indexed by elements $\sigma$ in $G$:
$$\zeta_S(\sigma,s) = \sum_{(\frak a,S) =1, \sigma_{\frak a} = \sigma} N(\frak a)^{-s}$$

Let $U_{\p}$ be the subgroup of the minus $S$ units $\epsilon$ in $K$ with $|\epsilon|_v = 1$ for all places $v$ not dividing $\p$. Stark conjectures that there is an element $\epsilon$ in $U_{\p}$ with the following properties. First, for all $\sigma$ in $G$
$$\log|\epsilon^{\sigma}|_{\frak q} = - ~W~\zeta'_S(\sigma,0).$$
where $W$ is the number of roots of unity in $K$. This identity, if true, determines the absolute value of $\epsilon$ at all places of $K$, and hence the $S$ unit $\epsilon$ up to a root of unity. Stark further conjectures that the $W^{th}$ root of $\epsilon$ lies in an abelian extension of $k$. When $k = \Q$ and $K = \Q(\mu_N)$, $\epsilon$ is a Jacobi sum and its $N^{th}$ root is a Gauss sum.

Since
$\zeta'_S(\sigma,0) = (\log N\frak p)\zeta_{S'}(\sigma,0)$
where $S' = S - \{\frak p\},$ the conjecture asserts that $\epsilon$ is a minus unit generating the ideal $\prod_{\sigma}\sigma^{-1}(\frak q)^{W ~\zeta_{S'}(\sigma,0)}$. We can get more information about this minus unit by considering the derivative of the $p$-adic partial zeta functions, where I conjecture that
\begin{equation} \label{e:gross}
 \log_p(N_{K_{\frak q}/\Q_p}\epsilon^{\sigma}) = - ~W~\zeta'_p(\sigma,0).
 \end{equation}
 \newpage

\section{Regulators in integral group rings}

I formulated the $p$-adic Stark conjecture and the rank one refinement in $1979$. I spent the next year testing these conjectures computationally, in some cases where the base field $k$ is real quadratic. Together with my student Leslie Federer, I checked that the analogous statement was true for Iwasawa modules \cite{FeG}. But I couldn't make any progress on the general proof, so I turned my attention to other topics, like Heegner points on modular curves.

I came back to the subject in the summer of $1985$, when I participated in a seminar at Harvard with Barry Mazur and John Tate on what we called ``refined conjectures". They were interested in refining the conjecture of Birch and Swinnerton-Dyer for the L-series of elliptic curves (cf. \cite{MT}), and I was interested in refining Stark's conjecture for Artin L-series \cite{Gr2}. What emerged in both cases was a general format, involving some calculus in the integral group ring of a finite abelian group. This formulation implied my earlier conjectures on $p$-adic L-functions, but was independent of the choice of a prime $p$.

I first treated the zeta function of a global field $k$, as David Hayes was attending the seminar and he wanted to understand the conjecture for function fields. Let $S$ be a non-empty set of places $v$ of $k$, which contains all the Archimedean places.
The zeta function $\zeta_S$ of $k$ is defined by the Euler product
$$\zeta_S(s) = \prod_{v \notin S}(1 - Nv^{-s})^{-1}$$in the half plane where the real part of $s$ is greater than $1$. It has a meromorphic continuation to the entire complex plane with a simple pole at $s=1$. It vanishes to order 
$$n = n(S) = \#S - 1$$
at the point $s = 0$, and has leading term $(-h_SR_S/w)$ there. Here $h_S$  is the order of the ideal class group of the Dedekind domain $A_S$ of $S$-integers in $k$, $w$ is the number of roots of unity in $k$, and $R_S$ is a generalization of Dirichlet's regulator on the units. Let $U_S = A_S^*$ be the group of $S$-units and $X_S$ the subgroup of elements of degree $0$ in the free abelian group on the places of $S$. Both $U_S$ and $X_S$ are finitely generated of rank $n(S)$. The torsion subgroup of $A_S^*$ is cyclic of order $w$ and the group $X_S$ is torsion free. The regulator $R(S)$ is the determinant of the homomorphism $\lambda_S:U_S \rightarrow X_S \otimes \R$
defined by
$$\lambda_S(\epsilon) = \sum_{v \in S} \log|\epsilon|_v ~v.$$
This determinant is taken on bases for the free parts, which are oriented so that the result is positive.

The only blemish in this formula is the presence of $w$ in the denominator of the leading term, which is a reflection of the fact that $\zeta_S(s)$ has a simple pole at $s = 1$. We can rectify the latter by introducing a finite set $T$ of places, which is non-empty and disjoint from $S$, and considering the zeta function modified at $T$
$$\zeta_{S,T}(s) = \prod_{v \in T}(1 - Nv^{1-s}) \zeta_S(s).$$
This still has order $n = n(S)$ at $s = 0$, but the leading term is now $-h_{S,T}R_{S,T}/w_T$, where $h_{S,T}$ is the order of a modified ideal class group with a trivialization of the projective modules at $T$ and $R_{S,T}$ is the determinant of $\lambda_S$ on the subgroup $U_{S,T}$ of $S$-units which are congruent to $1$ at all places in $T$. If we choose $T$ large enough so that this subgroup is torsion-free, then $w_T = 1$ and there are no denominators to deal with.

The basic idea of ``refining" this formula is to consider other homomorphisms $U_{S,T} \rightarrow X_S \otimes G$, where $G$ is a finite abelian group. (Barry and John were considering height pairings with values in $G$.) The values $\log|\epsilon|_v$ in $\lambda_S$ come from the local components of the ad\`elic norm $N: \A^*/k^*\prod_v A_v^* \rightarrow \R^*$.
Associated to a finite abelian extension $L$ of $k$, which is unramified outside of $S$ with Galois group $G$, we have the reciprocity homomorphism of global class field theory
$$F: \A^*/k^*\prod_{v \notin S}A_v^* \rightarrow G.$$
Using this, we can define an analogous homomorphism $U_{S,T} \rightarrow X_S \otimes G$ by
$$\lambda_F(\epsilon) = \sum_{v \in S} F_v(\epsilon) ~v.$$
Here $F_v$ is the restriction of $F$ to the subgroup $k_v^*$ of $\A^*$.

The $n \times n$ matrix of $\lambda_F$, with respect to our oriented bases of the free abelian groups $U_{S,T}$ and $X_S$, consists of $n^2$ elements of $G$. Let $\Z[G]$ be the integral group ring of $G$ and let $I$ be the augmentation ideal of elements of degree zero. The map from $G$ to $I$ defined by $g \rightarrow (g) - (1)$ gives an isomorphism of abelian groups $G \cong I/I^2$. (I learned this cohomological fact from John Tate, when we were windsurfing on the Mystic Lakes.) Hence we may consider the elements of the matrix as elements of $I/I^2$, and define its determinant (by the usual polynomial of degree $n$ in $n$ variables) as an element of $I^n/I^{n+1}$, which we call $\det(\lambda_F)$. (Barry and John initially defined their regulators in the group $\Sym^n(G)$, but I eventually persuaded them to work with the group $I^n/I^{n+1}$ instead.)

The refined conjecture relates this regulator to an element $\theta_{S,T}$ in $\Z[G]$ coming from the abelian extension $L$. This element has the property that for all characters $\chi$ of $G$
$$\chi(\theta_{S,T}) = L_{S,T}(\chi^{-1}, 0).$$
It lies in the integral group ring by results of Deligne-Ribet \cite{DR} and Cassou-Nogues \cite{C}. I conjectured that $\theta_{S,T}$ lies in the ideal $I^n$ of $\Z[G]$, with $n = \#S - 1$ and that modulo $I^{n+1}$ it is congruent to $-h_{S,T} \det(\lambda_F)$.

Henceforth we will assume that $n \geq 1$ as it is an easy matter to check the conjecture in the few cases when $n = 0$. Then $I^n/I^{n+1}$ is a finite abelian group, with the same prime divisors as $G$. If $G$ is cyclic, $I^n/I^{n+1}$ is isomorphic to $G$ for all $n \geq 1$. In general, the group $I^n/I^{n+1}$ stabilizes for large $n$, depending on $G$. For example, if $G = (\Z/p\Z)^2$ then $I^n/I^{n+1} = (\Z/p\Z)^{n+1}$ for $n \leq p$, and $I^n/I^{n+1} = (\Z/p\Z)^{p+1}$ for $n \geq p$.

In the number field case, this conjecture reduces to a question on $2$-groups. Indeed, if $k$ has a complex place $v$, then $F_v(\epsilon) = 1$ in $G$ for all $S$-units $\epsilon$. Hence an entire row of the matrix of $\lambda_F$ lies in $I^2$ and the determinant lies in $I^{n+1}$. In this case, we also have $L_{S,T}(\chi,0) = 0$ for all $\chi$, so $\theta_{S,T} = 0$  in $\Z[G]$ and the conjecture is true. If we assume that $k$ is totally real, the conjecture is also true for trivial reasons unless $L$ is totally complex. Then for each real place $v$, $F_v(\epsilon)$ has order $1$ or $2$ in $G$, depending on the sign of the unit $\epsilon$ at $v$. Hence $2\det(\lambda_F)$ lies in $I^{n+1}$. Similarly, $L_{S,T}(\chi, 0) = 0$ unless $\chi$ is totally odd. I was able to check the conjecture in many of these cases \cite[\S5]{Gr2}, using some subtle results of Deligne and Ribet on the $2$-divisibility of L-values \cite{DR}, and Hirose later proved it in general \cite{H}.

\section{From integral group rings to $p$-adic L-functions} \label{s:integral}

As mentioned, the formulation in the previous section does not imply the $p$-adic Stark conjecture, owing to Archimedean considerations. But we can modify it slightly to do so \cite[\S8]{Gr2}. Here we treat the case when $k$ is the totally real subfield of the CM field $K$, and $\chi$ is the non-trivial character of the Galois group. Let $L$ be a finite abelian extension of $k$ with group $G$. We assume that $L$ does not contain $K$, and let $S$ be a finite non-empty set of places of $k$ which contains all Archimedean places and all places ramified in $LK$. Finally, let $T$ be a finite set of places disjoint from $S$ which separates the roots of unity in any abelian extension of $k$ which is unramified outside of $S$.

We consider the one dimensional torus $H$ over $k$ which is split by $K$. This can be defined over the ring $A_S$ of $S$ integers, and the abelian group $H(A_S)$ is finitely generated, of rank equal to the number $n$ of places $\frak p$ in $S$ which split in $K$. This is the same rank as the free group $X_{S,K}^-$. Note that all places $\frak p$ which split are finite, as $K$ is a totally complex extension of the totally real field $k$. The Galois group $G$ of $L/k$ is also the Galois group of the abelian extension $LK/K$, so by global class field theory we have a surjective homomorphism $F: \A_K^*/K^*\prod_{v \notin S}A_{K,v}^* \rightarrow G$. This gives a homomorphism
$$\lambda_F: H(A_S)_T \rightarrow X_{S,K}^-\otimes G$$
defined by $\lambda_F(\epsilon) = \sum_{v \in S} F_v(\epsilon) ~v.$ Choosing bases for the free groups of rank $n$, we obtain an $n \times n$ matrix with entries in $G = I/I^2$. Hence $\det (\lambda_F)$ lies in $I^n/I^{n+1}$.

The conjecture for the torus $H$ relates this regulator with a theta element $\theta_{S,T}$ in $\Z[G]$ made from the values of abelian L-functions. For every character $\rho$ of $G$, we insist that $\rho(\theta_{S,T}) = L_{S,T}(\chi.\rho^{-1},0)$. Again, it is a non-trivial theorem that this element is integral. The refined conjecture for the torus $H$ states that this element lies in the ideal $I^n$, and that modulo $I^{n+1}$ is the product of $\det(\lambda_F)$ by an integer coming from the complex L-function of $\chi$.

If we assume that the set $S$ contains all prime factors of $p$, we can take $L = L_m$ as the $m$th cyclotomic extension of $k$ with Galois group isomorphic to the cyclic group $\Z/p^m\Z$. The refined conjecture is compatible with field extension, so we get a sequence of elements $\theta_{S,T}$ in $I^n/I^{n+1} = \Z/p^m/\Z$ converging to a an integral measure in $\Z[[\Z_p]]$. Its integral against the cyclotomic character gives the $p$-adic L-function $L_p(\chi\omega,s)$, and our group ring conjecture reduces to the $p$-adic Stark conjecture.

At this point, I was happy with the formulation of the $p$-adic Stark conjecture and its refinements, but I didn't see how to make any progress on it. So I turned it over to my talented students. In his PhD thesis, Henri Darmon showed that the methods which Kolyvagin had introduced to construct classes and annihilators in the Selmer group worked perfectly in the setting of the refined conjectures. And in his undergraduate thesis, Samit Dasgupta began his study of the original conjecture of Stark, as well as its $p$-adic analog. I'll let him continue with the story.

\section{Explicit formula} \label{s:exact} 

In the spring of 1998 I (Samit Dasgupta) asked Dick Gross, my undergraduate advisor, to suggest a topic and an advisor for a senior honors thesis.  He suggested that I read Tate's book on Stark's conjectures \cite{T} and work with Brian Conrad.  Stark's conjectures captivated me, and after writing my thesis \cite{dasthesis} I continued to ponder them while in graduate school.  In the summer of 2001, I attended a brilliant series of lectures by Henri Darmon \cite{darmonbook} on his recent conjectures on ``Stark--Heegner points'' \cite{darmon}.  At the conference I discussed possible thesis topics with him, and he suggested that his constructions could potentially be generalized to the context of Stark units.  

My resulting Ph.D. thesis (supervised by both Darmon and Ken Ribet) considered the setting of a real quadratic field $k$ and a prime $p$  inert in $k$.  I stated a conjectural exact formula for the Brumer--Gross--Stark unit $u_p$ in $k_p^*$, generalizing Gross's conjecture (\ref{e:gross}) for $\N_{k_p/\Q_p}(u_p) \in \Q_p^*$.  One appeal of this conjecture is that I proved that the set of all Brumer--Stark units (for a fixed ground field $F$ and all primes $p$) together with some simple square roots generates the maximal abelian extension of $F$.  Thus my conjecture gives a $p$-adic viewpoint on Hilbert's 12th problem.

The central new perspective employed in my thesis is that of Greenberg and Stevens in their study of $\sL$-invariants and the conjecture of Mazur, Tate, and Teitelbaum \cite{GS}.  A key idea of theirs is that the measures on $\Z_p^*$ used to define $p$-adic $L$-functions can be extended to $\PP^1(\Q_p)$ and then lifted to 
\[ (\Z_p^2)'  = \Z_p^2 - p\Z_p^2, \]
which is a $\Z_p^*$-bundle over $\PP^1(\Q_p)$.
Such lifts are determined by $p$-adic families of modular forms called {\em Hida families}.  In particular, they discovered that the leading terms of $p$-adic $L$-functions at trivial zeroes are controlled by {\em infinitesimal deformations} of modular forms.

\section{The rank one Gross--Stark conjecture} \label{s:ddp}

In the fall of 2003, while finishing the writing of my thesis and applying for postdoctoral positions, I visited Princeton along with Darmon for a special semester on the Birch--Swinnerton-Dyer conjecture.  Darmon, Wiles, and Skinner gave a course in which they rotated giving lectures on modular forms, Galois representations, and arithmetic applications.  Wiles spoke on Fermat's Last Theorem, Darmon on Stark--Heegner points, and Skinner spoke on his work with Urban on the main conjecture of Iwasawa theory for elliptic curves.  The Skinner--Urban approach generalized the work of Mazur--Wiles and Wiles on the Iwasawa main conjecture for Deligne--Ribet $p$-adic $L$-functions.  Wiles had proved the main conjecture for totally real fields by applying Ribet's method in the context of Hida families of Hilbert modular forms.

Ribet established his method in the transformative paper \cite{ribet} in which he proved the converse to Herbrand's Theorem.  This result states that if $p$ is an odd prime and $k > 0$ is an even integer such that $k \not\equiv 0 \pmod{p-1}$, then $p \mid \zeta(1 - k)$ implies that the component of $\text{Cl}(\Q(\mu_p)) \otimes \F_p$ on which Galois acts via $\chi^{1-k}$ is nontrivial, where $\chi$ is the mod $p$ cyclotomic character.  Ribet's technique, following in the footsteps of Siegel, is to first recognize $\zeta(1-k)$ as the constant term of the Eisenstein series $E_k$ on $\SL_2(\Z)$.  Ribet shows that the divisibility $p \mid \zeta(1 - k)$ implies that $p$ is congruent modulo $p$ to a cuspidal eigenform $f$ of weight $2$ and nebentypus $\omega^{k-2}$.  To this eigenform is associated a Galois representation (at the time of Ribet's writing, the existence and properties of these Galois representations were not yet fully established, and Ribet used some delicate algebraic geometry to construct the representations with the properties he needed).  The congruence between the cusp form and the Eisenstein series implies that this Galois representation is residually reducible, and hence gives rises to a Galois cohomology class for the module $\F_p(\chi^{1-k})$.  Furthermore, local considerations concerning the Galois representation allow one to prove that this cohomology class is everywhere unramified.  Global class field theory relates the existence of such a nonzero cohomology class to the nonvanishing of the desired component of $\text{Cl}(\Q(\mu_p)) \otimes \F_p$.  Mazur--Wiles (and later Wiles) applied Ribet's technique in the context of $p$-adic families of modular forms in order to relate $p$-adic $L$-functions to the characteristic ideals of certain Galois modules as predicted by the Iwasawa main conjecture.

Given the thesis work I had just finished under Darmon and Ribet, and the lectures we were now hearing from Skinner, it was impossible for Darmon and I {\em not} to wonder if the Ribet/Mazur--Wiles approach to $p$-adic $L$-functions could be used to prove Gross's conjecture.   
If the $p$-adic $L$-function is divisible by $T^r$, where $T=0$ corresponds to evaluation at some $s \in \Z_p$, then Wiles proved that there is a family of cusp forms congruent to an Eisenstein series modulo $T^r$.  The Greenberg--Stevens approach suggests that at a trivial zero, one should study whether there is such a congruence modulo $T^{r+1}$.  This is indeed what we found; for simplicity, let us describe this over the ground field $\Q$.

Let $\chi$ be an odd Dirichlet character such that $\chi(p)=1$.  For  each odd $k \ge 1$, there is a weight $k$ Eisenstein series
\[ E_k(1, \chi)(z) = \frac{L(\chi, 1-k)}{2} + \sum_{n=1}^{\infty} \left(\sum_{d \mid n} \chi(d) d^{k-1} \right) q^n. \]
The ordinary $p$-stabilizations of these  forms 
\[ E_k^*(1, \chi)(z) = \frac{L_p(\chi \omega, 1-k)}{2} + \sum_{n=1}^{\infty} \left(\sum_{d \mid n, p\nmid d} \chi(d) d^{k-1} \right) q^n \]
interpolate into a $p$-adic family.
We next apply an observation of Serre.  Consider the similar family defined with $\chi$ replaced by the trivial character (so its specializations are classical for {\em even} $k \ge 2$):
\[ E_k^*(1, 1)(z) = \frac{\zeta_p(1-k)}{2} + \sum_{n=1}^{\infty} \left(\sum_{d \mid n, p\nmid d} d^{k-1} \right) q^n. \]
Define the normalization $G_k^* = \frac{2}{\zeta_p(1-k)} E_k^*(1,1)$ whose constant term is 1:
\[ G_k^* = 1 + \frac{2}{\zeta_p(1-k)} \sum_{n=1}^{\infty} \left(\sum_{d \mid n, p\nmid d} d^{k-1} \right) q^n. \]
The $p$-adic zeta function has a pole at $s=1$, which implies that the specialization of the family $G_k^*$ at $k=0$ is simply the constant $1$. This was Serre's observation.

We use this to define the linear combination:
\begin{equation} \label{e:fkdef}
 F_k^* = E_k^*(1, \chi) - \frac{L_p(\chi\omega, 1-k)}{L(\chi, 0)} \cdot E_1(1, \chi) G_k^*. 
 \end{equation}
We make some observations: 
\begin{enumerate}
\item The coefficients have been chosen to cancel constant terms.  The form $F_k^*$ is cuspidal at infinity.
\item Since $L_p(\chi\omega, 1-k)$ has a trivial zero at $k=1$, the form $F_k^*$ specializes to $E_1^*(1, \chi)$ in weight 1.
\item In fact, since $G_0^* = 1$, the family $F_k^*$ is a Hecke eigenform in an infinitesimal neighborhood of weight $k=1$, i.e. it remains an eigenform if we expand as a Taylor series in $(k-1)$ and reduce modulo $(k-1)^2$.
\end{enumerate} 
The ``modulo $(k-1)^2$ Hecke eigenvalues'' of $F_k^*$ can be computed explicitly.  For prime $\ell \nmid p N$, where $N$ is the conductor of $\chi$, the eigenvalue is the same as that for $E_k^*(1,  \chi)$, namely 
\[ 1 + \chi(\ell) \ell^{k-1} \equiv 1 + \chi(\ell) + \chi(\ell) \log_p(\ell)(k - 1) \pmod{(k-1)^2}. \]
For the $U_p$-operator, the eigenvalue is more interesting:
\[ 1 - \frac{L_p'(\chi\omega,0)}{L(\chi, 0)} (k-1) \pmod{(k-1)^2}. \]
The ratio  $- L_p'(\chi\omega,0)/L(\chi, 0)$ is exactly that which appears in Gross's conjecture.
The method of Greenberg--Stevens shows how to use the Galois representation associated to a cuspidal family 
to relate the derivative of the $U_p$-operator of the family to a certain Galois cohomology class.
In our setting, the Galois cohomology class can be related to the Brumer--Gross--Stark unit $u_\chi$ via Kummer theory.  Therefore, all the pieces are in place to prove Gross's conjecture.  However there is a problem---the family $F_k^*$ is cuspidal at infinity, but it is not  cuspidal.

I was at Harvard as a postdoc at this point, and I began discussing with Robert Pollack at Boston University about how to modify the construction to make a cuspidal family.  Wiles's method is to apply a  certain Hecke operator, to pass from ``cuspidal at infinity'' to ``cuspidal'', but this operator actually annihilates our form.  We discovered the solution after Pollack wrote on the blackboard
\begin{quote}
``What would Ribet do?''
\end{quote}
I went home and studied Ribet's method of proving the cuspidality of the form in his original construction.  In Ribet's situation, it is easy to enumerate the Eisenstein series on $\Gamma_0(p)$ and to eliminate the possibility that the form constructed is any of these using congruences.  I then realized the problem in our situation---there is another Eisenstein family that we {\em cannot} eliminate using congruences.  Namely, the family 
\[ E_k^*(\chi, 1) = \sum_{n=1}^{\infty} \left(\sum_{d \mid n, p \nmid d} \chi(n/d)d^{k-1} \right)q^n \]
also specializes to $E_1^*(1, \chi)$ in weight $1$, and hence an appropriate multiple must be added to the linear combination (\ref{e:fkdef}).  Define
\begin{equation} \label{e:hkdef}
 H_k^* = E_k^*(1, \chi) - \frac{L_p(\chi\omega, 1-k)}{L(\chi, 0)} \cdot E_1(1, \chi) G_k^* - \frac{L_p(\chi\omega, 1-k)L(\chi^{-1},0)}{L_p(\chi^{-1}\omega, 1-k)L(\chi,0)}E_k^*(\chi, 1). 
 \end{equation}
It turns out that this linear combination has vanishing constant terms at both $\infty$ and $1/p$.  If
\begin{equation} \label{e:linv}
 \frac{L_p'(\chi\omega, 0)}{L(\chi, 0)} + \frac{L_p'(\chi^{-1}\omega, 0)}{L(\chi, 0)} \neq 0,
 \end{equation}
then one can define a Hecke operator $t$ such that $t\cdot H_k^*$ is a cuspidal ordinary family of modular forms that is normalized in the sense that $a_1(H_k^*) = 1$.  We also retain the property that $H_k^*$ is a ``mod $(k-1)^2$-eigenform'', though with slightly different Hecke eigenvalues than listed for $F_k^*$ above.
This whole construction works more generally in the context of Hilbert modular forms.  Combining this construction with ideas from Greenberg--Stevens, we arrived at the following result (\cite[Theorem 2]{ddp}). 

\begin{theorem} \label{t:ddp} Let $F$ be a totally real field, $\p$ a prime ideal of $F$ above $p$, and $\chi$ a totally odd character of $F$ such that $\chi(\p)=1$.  Suppose that Leopoldt's conjecture holds for $F$ at $p$, and that condition (\ref{e:linv}) holds as well.  Let $K$ be a CM field containing the fixed field of $\chi$, and fix a prime factor $\q$ of $\p$ in $K$.
Let $u_\chi$ be a generator of the 1-dimensional $\bC_p$-vector space $(\O_K[1/\p]^* \otimes \bC_p)^{\chi^{-1}}$.  Then
\begin{equation} \label{e:chiform}
 \frac{\log_p(N_{K_{\q}/\Q_p}(u_\chi))}{\ord_\q(u_\chi)} = - \frac{L_p'(\chi\omega, 0)}{L(\chi, 0)}. \end{equation}
\end{theorem}
It is not difficult to see that (\ref{e:chiform}) is equivalent to Gross's Conjecture (\ref{e:gross}) stated in \S\ref{s:abelian}.  The assumption of Leopoldt's conjecture for $(F,p)$ ensures that the $p$-adic zeta function of $F$ has a pole at $s=1$ (the equivalence of these facts was proven by Colmez \cite{colmez}).  This was used in the construction of the family $G_k^*$ specializing to the constant 1 in weight $k=0$.

\section{The higher rank case}

In February 2012, I visited UCLA to give some talks and met Kevin Ventullo, a graduate student of Chandrashekhar Khare.  Kevin was interested in my work with Darmon and Pollack.  He had the nice idea that one could remove the assumption of Leopoldt's conjecture from Theorem~\ref{t:ddp} by constructing the  family $G_k^*$ via other means.  Using an algebro-geometric perspective on Hilbert modular forms, one defines a lift of the Hasse invariant in high enough weight, and then constructs a $p$-adic family by taking powers of this form.

I was happy that Kevin had removed the assumption of Leopoldt's conjecture from our result, and thought that his work would make for an excellent thesis. I was extremely surprised when, several months later, he sent me a manuscript that removed not only this assumption, but also the assumption (\ref{e:linv}) as well!  In other words, Kevin had made the proof of the rank 1 abelian Gross--Stark conjecture completely unconditional.

Kevin's idea to remove the condition (\ref{e:linv}) was to simply imagine what happens if the condition fails.  In this case, the form $H_k^*$ is {\em not} a normalized eigenform modulo $(k-1)^2$, which  we had used to  provide a homomorphism from Hida's Hecke algebra to the ring of dual numbers $\bC_p[T]/T^2$.  However, one can study the Hecke orbit of this form, and analyze the action of the Hecke algebra on this space.  This provides a homomorphism from the Hecke algebra to a matrix algebra.  Somewhat miraculously to me, the image of this homomorphism once again turns out to be isomorphic  to the ring of dual numbers!  The proof can therefore proceed as before.

Upon reading Kevin's proof, I realized that the same idea could be used to approach the higher rank Gross conjecture.  With the notation as in \S\ref{s:padicconj}, let $\chi$ be a (1-dimensional) character of the abelian group $G = \Gal(K/k)$.
 Gross conjectured that $L_p(\chi\omega, s)$ vanishes to order $r = r_p(\chi)$ at $s=0$, where $r$ is the number of primes $\mathfrak{p}$ above $p$ in $k$ such that $\chi(\mathfrak{p}) = 1$.  Note that by Wiles' proof of the Iwasawa main conjecture (see also \cite{pcsd} and \cite{spiess}, which provide an alternate approach that in particular handles the case $p=2$), one knows that $L_p(\chi\omega, s)$ vanishes to order {\em at least} $r$ at $s=0$.  As described in \S\ref{s:padicconj}, Gross furthermore conjectured an exact formula for $L_p^{(r)}(\chi\omega, 0)$. 

When $r > 1$, the strategy employed in \cite{ddp} and described in \S\ref{s:ddp} encounters difficulties.
It is not possible to write down a single form that is a Hecke eigenform in an infinitesimal neighborhood of weight 1.  However, it is possible to write down a collection of forms that are Hecke stable, and to write down a homomorphism from the Hecke algebra into a certain ring of matrices.  In this case, the image of the homomorphism turns out to be a generalization of the ring of dual numbers that contains an infinitesimal variable for each prime $\mathfrak{p}$.  See \cite{dkv} for the precise definition.  This homomorphism also ``sees'' the leading term of the $p$-adic $L$-function.  I was therefore quite hopeful that following the methods of \cite{ddp} one could prove the higher rank conjecture.

I contacted Kevin about this idea, and of course he had already discovered the exact same construction.  We decided to join forces.  Despite much effort, we became stuck on completing the algebraic half of the proof.  In \cite{ddp}, we had used the Galois representations associated to cuspidal Hilbert modular forms to construct a Galois cohomology class that is {\em cyclotomic} in the sense that its restriction to a decomposition group at $p$ is a multiple of the logarithm of the cyclotomic character.  Applying the same strategy directly in the rank $r$ case, one needs to construct $r$ linearly independent such classes.  However, in this case we were unable to construct {\em any} cyclotomic classes; the rank 1 setting provided a certain rigidity that was now lacking.  After many months of collaborative effort, our progress stalled, and I began to focus on other projects.

In March 2013, Mahesh Kakde visited me at UC Santa Cruz to discuss a separate project related to the factorization of $p$-adic $L$-series \cite{factor}.  As an aside, he mentioned that the number theory group in London had run a seminar studying my paper with Darmon and Pollack, and he was interested in any thoughts I had on the higher rank case.  I described my progress with Kevin and where we were stuck.  After he returned to London, we began to discuss the problem over email.  We made some important progress during a trip I made for the Iwasawa 2015 conference in London; we figured out a rank 2 setting under certain additional assumptions that allowed us to mimic the ``rigidity'' that holds in the rank 1 case.  

After returning to Santa Cruz, I gave a graduate course on the proof of Dasgupta--Darmon--Pollack.  I wanted to give a down-to-earth explanation of the connection between units and cohomology classes that avoided the use of Poitou--Tate duality as in the paper.
 I came up with an explicit proof of the orthogonality between units and cohomology classes using nothing more than the adelic statements of class field theory.  While giving the lecture, it struck me that the same idea can be used in the higher rank case.  Instead of requiring ourselves to construct cyclotomic cohomology classes, we can simply consider the cohomology classes that occur and write down the statement that they annihilate a basis of the group of $p$-units.  The $r$ equations that emerge can be viewed as stating that a certain matrix (of cohomology classes evaluated locally at units) has row sums that vanish.  Of course, this forces the vanishing of the determinant of that matrix; expanding this equation gives precisely the desired result that the determinant of the logarithms of the norms of the units equals the leading term of the $p$-adic $L$-function.  See \cite[\S1.4]{dkv} for more details on this argument.  We arrived at a proof of Gross's conjecture on the (expected) leading term of Deligne--Ribet $p$-adic $L$-functions at $s=0$, described in \S\ref{s:padicconj} above.
 
 \begin{theorem} \label{t:gross} We have:
\begin{equation} \label{e:grosshr}
 \frac{L_p^{(r)}(\chi\omega, 0)}{r! L(\chi, 0)} =  R_p(\chi) \prod_{\chi(\mathfrak{p}) \neq 1} (1 - \chi(\mathfrak{p})).
 \end{equation}
 \label{e:g2}
\end{theorem}
 Our resulting paper \cite{dkv} was one of the most interesting collaborations I have been involved with---while all three of us contributed significantly to the project, Mahesh and Kevin never actually met until after the completion of the paper!

 In (\ref{e:grosshr}),  $R_p(\chi)$ is Gross's regulator defined in (\ref{e:rpdef}).  The non-vanishing of this  quantity is known as the Gross--Kuz'min conjecture, and it is known to follow from standard conjectures in transcendence theory (such as the $p$-adic Schanuel conjecture).  In view of Theorem~\ref{t:gross}, the non-vanishing of $R_p(\chi)$ is equivalent to the statement that 
$L_p(\chi\omega,s)$ vanishes to order exactly $r$ at $s=0$.

\section{The Integral Conjectures} \label{s:ic}
 
 Ever since my work with Darmon and Pollack in the mid 2000s, I had tried to establish an integral version of the constructions in order to prove Gross's integral Stark conjecture described in \S\ref{s:integral}, at least in the rank 1 case.  Technical problems always prevented me from succeeding, in particular the proper construction of the analogue of the auxiliary form $G_k^*$ appearing in (\ref{e:fkdef}).  After Kevin's work, I realized that a more theoretical approach to the construction, avoiding Serre's observation using the pole of the $p$-adic zeta function, could be fruitful.  However, difficulties in enacting Kevin's technique integrally, in particular in the presence of general level structure, remained.
 
 In January of 2019 Mahesh visited me at Duke, where I had moved a few months before, and showed me an interesting paper of David Burns \cite{burns}.  This built upon Burns' previous joint work with Masato Kurihara and Takamichi Sano \cite{bks}.  Burns proved an equivariant main conjecture that implied the Brumer--Stark conjecture under the assumption that the $\mu$-invariant of $F$ vanishes and that the Gross--Kuz'min conjecture holds (i.e. the non-vanishing of Gross's regulator).  
 Mahesh proposed that perhaps we could remove the $\mu$-invariant condition from Burns' result. 
 
 In fact we tried to be a little more ambitious.  Rather than working Iwasawa theoretically, we considered {\em group-ring} modular forms for $F$.  Given a finite CM abelian extension $K/F$ with $G =\Gal(K/F)$, these are defined as follows.  We consider Hilbert modular forms $f$ over $F$ with level equal to the conductor of $K$ and weight equal to some odd integer $k \ge 1$, with Fourier coefficients in the group ring $\Z_p[G]$, such that for each odd character $\chi$ of $G$ the form $\chi(f)$ has nebentypus $\chi$.  The fact that the form $f$ has coefficients in $\Z_p[G]$ implies that the forms $\chi(f)$ satsify certain congruences modulo $p$.

An important ingredient in the work of Burns--Kurihara--Sano is a study of the modules living in Tate sequences originally defined by Ritter and Weiss \cite{RW}.  Mahesh and I are indebted to Sano for making us aware of the Ritter--Weiss papers and suggesting that they could be useful in our work.
Another timely occurence is that I welcomed a new postdoc to Duke, Jesse Silliman, who is an expert on the algebraic geometry of Hilbert modular varieties.  Jesse was able to solve my difficulties regarding the existence of group-ring valued modular forms generalizing the $G_k^*$ \cite{silliman}.  

The pieces were in place to solve the Brumer--Stark and integral Gross--Stark conjectures, but once again there was an important technical difficulty.  In applying Ribet's method, one has to show that the $p$-adic cohomology classes constructed are unramified at $p$.  By class field theory, the ramification is controlled by a class group, which is finite.  When working over a field of characteristic 0 such as $\Q_p$, tensoring with the field annihilates this finite ramification.   However, when working over $\Z_p$, this finite ramification is an important obstacle.  In Wiles' proof of the Iwasawa main conjecture, he uses an ingenious method of twisting by auxiliary characters of conductor prime to $p$.  But again this technique introduces finite error terms, and extra conditions are necessary to achieve a full unconditional result (see \cite{greither}).

The issue is that the cohomology class that one constructs by directly applying the analogue of the congruence proved by Wiles in his study of the main conjecture {\em actually is} ramified at $p$.  Mahesh and I discovered that fortunately there is an  ``extra congruence'', arising as a shadow of the trivial zeroes of $p$-adic $L$-functions at $s=0$.  More precisely, we defined a non-zerodivisor $x \in \Z_p[G]^-$ such that there is a congruence between a cusp form and an Eisenstein series modulo $x \theta$, where $\theta$ is the Stickelberger element, rather than just modulo $\theta$.  One then obtains a cohomology class in a larger module than if one only used the congruence modulo $\theta$.  This cohomology class is ramified at $p$, but we show that modding out by the image of the inertia groups at $p$ cuts down the size of the module (in the sense of Fitting ideals) by a factor of precisely this same factor $x$.  One therefore obtains the desired cohomology class that is unramified at $p$ and generates a module with Fitting ideal $(x \theta)/x = \theta$.

As a corollary of our results, we obtained the following (see \cite{dk} and \cite{dk2}):

\begin{theorem} We have the following:
\begin{itemize}
\item  The Brumer--Stark conjecture holds over $\Z[1/2]$.  
\item  With notation as in \S\ref{s:integral}, assume that the set $S$ contains exactly $n=1$ prime $\mathfrak{p}$ of $k$ that splits completely in $K$.  Suppose $\mathfrak{p}$ lies above $p$ and that $p$ is odd.  Then 
the integral Gross--Stark conjecture  (i.e.~the conjecture of \S\ref{s:integral}) holds in $I/I^2 \otimes \Z_p$.
\item  Again assume that the prime $\mathfrak{p}$ of $k$ splits completely in $K$.  Assume that $p \nmid 2 \disc(F)$.
The exact formula for Brumer--Stark units in $k_\mathfrak{p}^*$ given in \cite{dd} and \cite{dasshin}, and mentioned in \S\ref{s:exact} above, 
holds up to a root of unity.
\end{itemize}
\end{theorem}

As described above, the third bullet point yields an explicit $p$-adic construction of the maximal abelian construction of any totally real field.

We conclude by noting that very recently,  Bullach, Burns, Daoud, and Seo \cite{bullach} make the nice observation that my results with Kakde  on the Brumer--Stark conjecture (more precisely, the ``Strong Brumer--Stark conjecture'' concerning Fitting ideals conjectured by Kurihara) imply the minus part of the Equivariant Tamagawa Number Conjecture for $K/k$ over $\Z[1/2]$.  The ETNC is essentially the strongest result in this direction, and is known to imply most, if not all conjectures in this setting.  In particular, one now has a proof of Gross's integral conjecture in $I^n/I^{n+1} \otimes \Z[1/2]$.

What remains is handle the troublesome prime $p=2$ and obtain results over $\Z$ rather than $\Z[1/2]$. The main problem in this setting is that Ribet's method encounters serious difficulties when attempting to construct extensions of a representation $\rho_1$ by a representation $\rho_2$ when the representations are ``residually indistinguishable'', i.e. when their reductions modulo $p$ are isomorphic.  In the case of the minus side for $K/k$, we are considering totally odd characters $\chi$ and trying to construct extensions of 1 by $\chi$.  Since $\chi(c) = -1 \not\equiv 1 \pmod{p}$ when $p$ is odd, Ribet's method works well.  New ideas are necessary in the residually indistinguishable case; this is curently work in progress with Kakde, Jesse Silliman, and Jiuya Wang.

\def\noopsort#1{}
\providecommand{\bysame}{\leavevmode\hbox to3em{\hrulefill}\thinspace}

\end{document}